\documentstyle[12pt]{article}
\def \qed{\hfill $\vrule height 2.5mm  width 2.5mm depth 0mm $}
\def\beq{\begin{eqnarray}}
\def\eeq{\end{eqnarray}}
\def\nn{\nonumber}
\def\wt{\widetilde}

\def\ld{\lambda}

\def\al{\alpha}
\def\wh{\widehat}

\def\mod{\hbox{\rm mod~}}

\def\ds{\displaystyle}
  
\def\qed{\hfill $\vrule height 2.5mm  width 2.5mm depth 0mm $}

\def\neweq{\setcounter{equation}{0}}
\newtheorem{theorem}{Theorem}[section]

\newtheorem{cor}[theorem]{Corollary}
\newtheorem{de}[theorem]{Definition}

\newtheorem{rem}[theorem]{Remark}

\def\m@th{\mathsurround=0pt}

\def\Fsquare(#1,#2){
\hbox{\vrule$\hskip-0.4pt\vcenter to #1{\normalbaselines\m@th
\hrule\vfil\hbox to #1{\hfill$#2$\hfill}\vfil\hrule}$\hskip-0.4pt
\vrule}}

\def\Addsquare(#1,#2){\hbox{$
    \dimen1=#1 \advance\dimen1 by -0.8pt
    \vcenter to #1{\hrule height0.4pt depth0.0pt\vss%
    \hbox to #1{\hss{%
    \vbox to \dimen1{\vss%
    \hbox to \dimen1{\hss$~#2~$\hss}%
    \vss}\hss}%
    \vrule width0.4pt}\vss%
    \hrule height0.4pt depth0.0pt}$}}

\def\Htwobox(#1,#2){%
    \Fsquare(0.4cm,#1)\Addsquare(0.4cm,#2)}

\def\Hthreebox(#1,#2,#3){%
    \Fsquare(0.4cm,#1)\Addsquare(0.4cm,#2)\Addsquare(0.4cm,#3)}

\def\Hfourbox(#1,#2,#3,#4){%
    \Hthreebox(#1,#2,#3)\Addsquare(0.4cm,#4)}

\def\Hfivebox(#1,#2,#3,#4,#5){%
    \Hfourbox(#1,#2,#3,#4)\Addsquare(0.4cm,#5)}

\def\Vfivebox(#1,#2,#3,#4,#5){%
    \normalbaselines\m@th\offinterlineskip
    \vtop{\hbox{\Fsquare(0.4cm,#1)}
          \vskip-0.4pt
          \hbox{\Fsquare(0.4cm,#2)}
          \vskip-0.4pt
          \hbox{\Fsquare(0.4cm,#3)}
          \vskip-0.4pt
          \hbox{\Fsquare(0.4cm,#4)}
          \vskip-0.4pt
          \hbox{\Fsquare(0.4cm,#5)}}}

\def\Vfourbox(#1,#2,#3,#4){%
    \normalbaselines\m@th\offinterlineskip
    \vtop{\hbox{\Fsquare(0.4cm,#1)}
          \vskip-0.4pt
          \hbox{\Fsquare(0.4cm,#2)}
          \vskip-0.4pt
          \hbox{\Fsquare(0.4cm,#3)}
          \vskip-0.4pt
          \hbox{\Fsquare(0.4cm,#4)}}}

\def\Vthreebox(#1,#2,#3){%
    \normalbaselines\m@th\offinterlineskip
    \vtop{\hbox{\Fsquare(0.4cm,#1)}
          \vskip-0.4pt
          \hbox{\Fsquare(0.4cm,#2)}
          \vskip-0.4pt
          \hbox{\Fsquare(0.4cm,#3)}}}

\def\Vtwobox(#1,#2){%
    \normalbaselines\m@th\offinterlineskip
    \vtop{\hbox{\Fsquare(0.4cm,#1)}
          \vskip-0.4pt
          \hbox{\Fsquare(0.4cm,#2)}}}

\def\Twoone(#1,#2,#3){%
    \hbox{
    \normalbaselines\m@th\offinterlineskip
    \vtop{\hbox{\Htwobox({#1},{#2})}
          \vskip-0.4pt
          \hbox{\Fsquare(0.4cm,#3)}}}}

\def\Threeone(#1,#2,#3,#4){%
    \normalbaselines\m@th\offinterlineskip
    \vtop{\hbox{\Hthreebox({#1},{#2},{#3})}
          \vskip-0.4pt
          \hbox{\Fsquare(0.4cm,#4)}}}

\def\Fourone(#1,#2,#3,#4,#5){%
    \normalbaselines\m@th\offinterlineskip
    \vtop{\hbox{\Hfourbox({#1},{#2},{#3},{#4})}
          \vskip-0.4pt
          \hbox{\Fsquare(0.4cm,#5)}}}

\def\Threetwo(#1,#2,#3,#4,#5){%
    \normalbaselines\m@th\offinterlineskip
    \vtop{\hbox{\Hthreebox({#1},{#2},{#3})}
          \vskip-0.4pt
          \hbox{\Htwobox({#4},{#5})}}}

\def\Twotwo(#1,#2,#3,#4){%
    \normalbaselines\m@th\offinterlineskip
    \vtop{\hbox{\Htwobox({#1},{#2})}
          \vskip-0.4pt
          \hbox{\Htwobox({#3},{#4})}}}

\def\Twooneone(#1,#2,#3,#4){%
    \normalbaselines\m@th\offinterlineskip
    \vtop{\hbox{\Htwobox({#1},{#2})}
          \vskip-0.4pt
          \hbox{\Fsquare(0.4cm,#3)}
          \vskip-0.4pt
          \hbox{\Fsquare(0.4cm,#4)}}}

\def\Twooneoneone(#1,#2,#3,#4,#5){%
    \normalbaselines\m@th\offinterlineskip
    \vtop{\hbox{\Htwobox({#1},{#2})}
          \vskip-0.4pt
          \hbox{\Fsquare(0.4cm,#3)}
          \vskip-0.4pt
          \hbox{\Fsquare(0.4cm,#4)}
          \vskip-0.4pt
          \hbox{\Fsquare(0.4cm,#5)}}}

\def\Twotwoone(#1,#2,#3,#4,#5){%
    \normalbaselines\m@th\offinterlineskip
    \vtop{\hbox{\Htwobox({#1},{#2})}
          \vskip-0.4pt
          \hbox{\Htwobox({#3},{#4})}
              \vskip-0.4pt
          \hbox{\Fsquare(0.4cm,#5)}}}

\def\Threeoneone(#1,#2,#3,#4,#5){%
    \normalbaselines\m@th\offinterlineskip
    \vtop{\hbox{\Hthreebox({#1},{#2},{#3})}
          \vskip-0.4pt
          \hbox{\Fsquare(0.4cm,#4)}
              \vskip-0.4pt
          \hbox{\Fsquare(0.4cm,#5)}}}
\begin{document}

\begin{center}
{\textbf{BETHE'S STATES FOR GENERALIZED $XXX$\\ AND $XXZ$ MODELS}}
\end{center}

\begin{center}
{Anatol N.\ KIRILLOV \& Nadejda LISKOVA}
\end{center}
\begin{center}
{\small {\it Graduate School of Mathematics, Nagoya University}} \\
{\small {\it Chikusa--ku, Nagoya 486--8602, Japan}}\\ {\small
{\it and}}
 \\
{\small {\it Steklov Mathematical Institute,}} \\ {\small {\it
Fontanka 27, St.Petersburg, 191011, Russia}}
\end{center}

\vskip -1cm

\begin{abstract}For any rational number $p_0\ge 1$ we prove an
identity of Rogers--Ramanujan--Gordon--Andrews' type. Bijection
between the space of states for $XXZ$ model and that of $XXX$
model is constructed.
\end{abstract}

\section{Introduction}
\label{intro}
\neweq

The main goal of our paper is to study a combinatorial
relationship between the space of states for generalized $XXZ$
model and that for $XXX$ one. In our previous paper \cite{KL} we
gave a combinatorial description of states for generalized $XXZ$
model in terms of the so--called $sl(2)$--$XXZ$ rigged
configurations. On the other hand it is well--known that when the
anisotropy parameter $p_0$ of $XXZ$ model goes to infinity then
the $XXZ$ model under consideration transforms to the $XXX$ one.
We are going to describe this transformation from combinatorial
point of view in the case when $p_0$ is an integer.

A combinatorial completeness of Bethe's states for generalized
$XXX$--model had been proved in \cite{K1} and appeared to be a
starting point for numerous applications to combinatorics of Young
tableaux and representation theory of symmetric and general
linear groups, see e.g. \cite{K2}. Here we mention only a
"fermionic" formula for the Kostka--Foulkes polynomials, see e.g.
\cite{K2}, and the relationship of the latter with
$\wh{sl}(2)$--branching functions $b_{\ld}^{k\Lambda_0}(q)$, see
e.g. \cite{K3}. We will show in Section~1, Theorem~\ref{t2.3} and
Remark~\ref{r2.5}, that $q$--counting of the number of $XXZ$
states using Bethe's ansatz approach \cite{TS,KR}, gives rise to
the Rogers--Ramanujan--Gordon--Andrews' type identity for any
rational number $p_0\ge 1$.

It seems an interesting problem to find a polynomial version of
the Rogers--Ramanujan type identity (\ref{2.12}) from our
Theorem~\ref{t2.3}.

Another question which we are interested in is to understand a
combinatorial nature of the limit
$$XXZ ~_{\buildrel\longrightarrow\over{p_0\to +\infty}}~ XXX.
$$

In Section~3 we describe a combinatorial rule which shows how the
$XXZ$--con\-fi\-gura\-tions fall to the $XXX$ pieces. For
simplicity we consider in our paper only the case
$p_0>\sum_ms_m$. General case will be considered elsewhere.

\section{Rogers--Ramanujan's type identity}
\label{rrti}
\neweq

 This paper is a continuation of our previous work \cite{KL}.
 Let us remind the main definitions, notation and results from \cite{KL}.

 For fixed $p_0\in{\bf R}$, $p_0\ge 1$ let us define (cf. \cite{TS})
 a sequence of real
 numbers $p_i$ and sequences of integer numbers $\nu_i,m_i,y_i,z_i$:
 \begin{eqnarray}
&&p_0:=p_0,~p_1=1,~\nu_i=\left[\frac{p_i}{p_{i+1}}\right],~
p_{i+1}=p_{i-1}-\nu_{i-1}p_i,~i=1,2,\ldots \\
&&y_{-1}=0,~y_0=1,~~y_1=\nu_0,~y_{i+1}=y_{i-1}+\nu_iy_i,~i=0,1,2,\ldots
\\
&&z_{-1}=0,~z_0=1,~z_1=\nu_1,~z_{i+1}=z_{i-1}+\nu_{i+1}z_i,~i=0,1,2,\ldots
\\
&&m_0=0,~m_1=\nu_0,~m_{i+1}=m_i+\nu_i,~i=0,1,2,\ldots \\
&& r(j)=i,~\hbox{\rm if}~m_i\le j<m_{i+1},~j=0,1,2,\ldots
 \end{eqnarray}
 It is clear that integer numbers $\nu_i$ define the decomposition of
 $p_0$ into continuous fraction
 $$p_0=[\nu_0,\nu_1,\nu_2,\ldots ]=\nu_0+\frac{1}{\nu_1+\ds\frac{1}{\nu_2+
 \ldots}}.
 $$
 Let us define (see Fig. 1) a piecewise linear function $n_j$, $j\ge 0$,
 \begin{equation}\label{2.5}
n_j:=y_{i-1}+(j-m_i)y_i,~~\hbox{\rm if}~~m_i\le j<m_{i+1}.
 \end{equation}

It is clear that for any integer $n>1$ there exists the unique
rational
 number $t$ such that $n=n_t$.

Let us introduce additionally the following functions (see
\cite{KL})
 \begin{equation}\label{2.6}
 q_j=(-1)^i(p_i-(j-m_i)p_{i+1}),~~\hbox{\rm if}~~m_i
\le j<m_{i+1},
 \end{equation}
 $$\Phi_{k,2s}=\left\{\begin{array}{ll}
 \ds\frac{1}{2p_0}(q_k-q_kn_{\chi}), &\mbox{if $n_k>2s$,}\\
 \ds\frac{1}{2p_0}(q_k-q_{\chi}n_k)+\ds\frac{(-1)^{r(k)-1}}{2}, &\mbox{if
 $n_k\le 2s$},
 \end{array}\right.
 $$
 where $2s=n_{\chi}-1$.

In what follows we assume that the anisotropy parameter $p_0\ge 1$
is a rational number, and $p_0=[\nu_0,\nu_1,\ldots ,\nu_{\alpha}]$
denotes  its decomposition into continuous fraction. It is not
difficult to see that
\begin{eqnarray*}
p_0&=&[\nu_0,\nu_1,\ldots ,\nu_{\al}]=\ds\frac{y_{\al
+1}}{z_{\al}},\\
\bar p_0&:=&[\nu_0,\nu_1,\ldots ,\nu_{\al
-1}]=\ds\frac{y_{\al}}{z_{\al -1}},
\end{eqnarray*}
where the numbers $\{ y_j\}_{j=0}^{\al +1}$ and
$\{z_j\}_{j=0}^{\al}$ are defined by (2.2) and (2.3)
correspondingly. We assume that if $\al >0$, then $\nu_{\al}\ge
2$.

\vskip 01.8cm \setlength{\unitlength}{0.4cm}
\begin{picture}(29,23)(-6.2,-1)
\put(0,0){\vector(0,1){24}} \put(0,0){\vector(1,0){25}}
\multiput(6,0)(0,1){3}{\line(0,1){0.5}}
\multiput(0,2)(1,0){6}{\line(1,0){0.5}}
\multiput(0,4)(1,0){8}{\line(1,0){0.5}}
\multiput(0,3)(1,0){12}{\line(1,0){0.5}}
\multiput(0,6)(1,0){10}{\line(1,0){0.5}}%
\multiput(0,8)(1,0){21}{\line(1,0){0.5}}
\multiput(0,13)(1,0){17}{\line(1,0){0.5}}
\multiput(0,17)(1,0){19}{\line(1,0){0.5}}
\multiput(0,21)(1,0){26}{\line(1,0){0.5}}
\multiput(8,0)(0,1){4}{\line(0,1){0.5}}%
\multiput(10,0)(0,1){6}{\line(0,1){0.5}}%
\multiput(12,0)(0,1){8}{\line(0,1){0.5}}
\multiput(17,0)(0,1){13}{\line(0,1){0.5}}
\multiput(19,0)(0,1){17}{\line(0,1){0.5}}
\multiput(21,0)(0,1){21}{\line(0,1){0.5}}
\multiput(6.5,3.25)(0.5,0.25){4}{\circle*{0.15}}
\multiput(12.5,8.5)(0.5,0.5){10}{\circle*{0.15}} \thicklines
\thicklines \put(0.8,0.4){\vector(2,1){5.2}}
\put(6,2){\vector(1,1){6}}
\put(12,3){\vector(1,2){9}}
\put(21,8){\line(1,4){4}} \put(-2,2){$y_{i-2}$}
\put(-2,3){$y_{i-1}$} \put(-4,4){$y_{i-1}\!\! +\!\! y_{i-2}$}
\put(-3.5,6){$y_i\! -\! y_{i-1}$} \put(-2,8){$y_{i}$}
\put(-3.5,13){$y_i\! +\! y_{i-1}$} \put(-3.5,17){$y_{i+1}\!\!
-\!\! y_i$} \put(-2,21){$y_{i+1}$} \put(4.5,-1.2){$m_{i-1}$}
\put(6.7,-1.2){$m_{i-1}\!\! +\!\! 1$} \put(9.9,-1.2){$m_i\!\!
-\!\! 1 $} \put(12.2,-1.2){$m_{i}$} \put(15.6,-1.2){$m_i\!\!
+\!\! 1$} \put(21.2,-1.2){$m_{i+1}$} \put(17.9,-1.2){$m_{i+1}\!\!
-\!\! 1$} \put(25,-1.2){$j$} \put(-1.5,23.5){$n_j$}
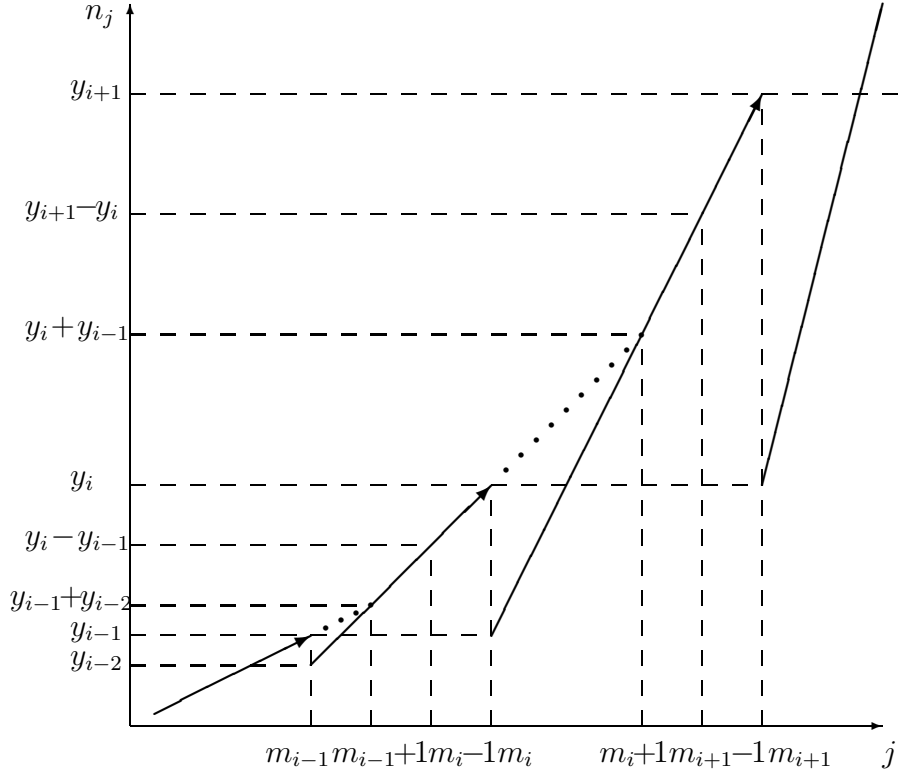
\put(-3,-3.5){Fig.1.~~Image of piece--wise linear function $n_j$
in the interval $[m_{i-1},m_{i+1}]$}
\end{picture}
\vskip 2cm

In order to formulate our main result of the paper \cite{KL}
about the number of Bethe's states for generalized $XXZ$ model,
let us consider the following symmetric matrix\break 
$\Theta^{-1}=(c_{ij})_{1\le i,j \le m_{\alpha +1}}$:
 \beq
&&i)~~c_{ij}=c_{ji}~~\hbox{\rm and}~~c_{ij}=0,
~~\hbox{\rm if}~~|i-j|\ge 2.\nn\\~\nn \\
&&ii)~~c_{j-1,j}=(-1)^{i-1},~~\hbox{\rm if}~~m_i\le j<m_{i+1}.\nn\\~\nn \\
&&iii)~~c_{jj}=\left\{\begin{array}{ll}
 2(-1)^i,&\mbox{if $m_i\le j<m_{i+1}-1,~~i\le \alpha$},\\
 (-1)^i, &\mbox{if $j=m_{i+1}-1$,~~$i\le\alpha$},\\
 (-1)^{\al +1},& \mbox{if $j=m_{\alpha +1}$}.\nn
 \end{array}\right.
 \eeq
\medbreak
\hskip -0.65cm{\bf Example 1} \ Let us take $p_0=\ds\frac{16}{7}$, 
then $p_0=[2,3,2]$, $\al =2$, $\bar p=[2,3]=\ds\frac{7}{3}$, 
\begin{eqnarray*}
&&\nu_0=2,~~\nu_1=3,~~\nu_2=2;\\
&&m_0=0,~~m_1=2,~~m_2=5,~~m_3=7;\\
&&y_0=1,~~y_1=2,~~y_2=7,~~y_3=16;\\
&&z_0=1,~~z_1=3,~~z_2=7.
\end{eqnarray*}
Therefore, $p_0=\ds\frac{y_{\al +1}}{z_{\al}}$, $\bar
p_0=\ds\frac{y_{\al}}{z_{\al -1}}$, and
$$n_j=\cases{j,& if $0\le j<2$,\cr
1+2(j-2), & if $2\le j<5$,\cr
2+7(j-5),& if $5\le j<7$,\cr
7+16(j-7), & if $7\le j$.\cr}~~~~~~~~~~~~~~~~~~~~~~~~~~~~~
$$
Finally,
 \medbreak
 $$\Theta =\frac{1}{16}\left(\begin{array}{rrrrrrr}
 9&7&5&3&2&1&1\\ 7&-7&-5&-3&-2&-1&-1\\ 5&-5&-15&-9&-6&-3&-3\\
 3&-3&-9&-15&-10&-5&-5\\ 2&-2&-6&-10&4&2&2\\ 1&-1&-3&-5&2&9&9\\
 1&-1&-3&-5&2&9&-7\end{array}\right),
 $$
 \vskip 0.8cm
 $$\Theta^{-1}=\left(\begin{array}{rrrrrrr}
1 & 1\\ 1 & -2 & 1\\ & 1 & -2 & 1\\
&& 1 & -1 &-1\\
  &&&-1&2&-1\\ &&&&-1&1&1\\ &&&&&1&-1
  \end{array}\right),
$$
and ~~$\det |\Theta^{-1}|=16$.

\vskip 0.5cm\hskip -0.65cm Note, see \cite{KL}, Theorem~4.7, that
the absolute value of the determinant $\det(\Theta^{-1})$ is
equal to $y_{\alpha +1}$, i.e. the numerator of $p_0$. \vskip
0.5cm

To continue, let us consider the matrix $E=(e_{jk})_{1\le j,k\le
m_{\alpha +1}}$, where
$$e_{jk}=(-1)^{r(k)}(\delta_{j,k}-\delta_{j,m_{\al +1}-1}\cdot\delta_{k,m_{\al
+1}}+\delta_{j,m_{\al +1}}\cdot\delta_{k,m_{\al +1}-1}).
$$
One can check that the vacancy numbers $P_j(\ld )$, see \cite{KL},
(3.9), can be computed as follows
$$P_j(\ld )+\ld_j=((E-2\Theta ){\wt\ld^t}+b^t)_j,~~1\le j\le m_{\al +1},\nn
$$
where the $j$-th component of vector $b=(b_1,\ldots ,b_{m_{\al
+1}})$ is defined by the following formula
$$b_j=(-1)^{r(j)}\left(n_j\left\{{\sum 2s_mN_m-2l\over p_0}\right\}
-\sum_m2\Phi_{j,2s_m}\cdot N_m\right),\nn
$$
and for any sequence of integer numbers $\ld =(\ld_1,\ldots
,\ld_{m_{\al +1}})$ we denote by $\wt\ld$ the sequence
$(\wt\ld_1,\ldots ,\wt\ld_{m_{\al +1}})$, where
$\wt\ld_j=(-1)^{r(j)}\ld_j$, $1\le j\le m_{\al +1}$.
 \medbreak
\begin{theorem} {\rm (\cite{KL})} The number of Bethe's states
$Z^{XXZ}(N,s~|~l)$ of the generalized $XXZ$ model is equal to
 \begin{equation}\label{2.7}
\sum_{\ld}\prod_j\left(
 \begin{array}{c}
    ((E-B){\wt\ld}^t+b^t)_j  \\
    \ld_j
 \end{array}\right),
 \end{equation}
where summation is taken over all sequences of non--negative
integer numbers $\ld =\{\ld_k\}_{k=1}^{m_{\al +1}}$ such
 that
 \begin{eqnarray*}
     &&\sum_{k=1}^{m_{\al +1}}n_k\ld_k=l,~~\ld_k\ge 0;  \\\nn \\
     &&{\wt\ld}=({\wt\ld}_1,\ldots .{\wt\ld}_{m_{\al
     +1}}),~~{\wt\ld}_j=(-1)^{r(j)}\ld_j,~~B=2\Theta.
 \end{eqnarray*}
 \end{theorem}

Recall, see \cite{KL}, or Section~3, that $N$ and $s$ in the
above formula for the number of states $Z^{XXZ}(N,s|l)$ denote
vectors $N=(N_1,\ldots ,N_k)$ and $s=(\underbrace{s_1,\ldots
,s_1}_{N_1},\ldots ,\underbrace{s_k,\ldots ,s_k}_{N_k})$, i.e.
$N_m$ is equal to the number of spins in the $XXZ$--chain which
are equal to $s_m$.

One of the main goal of the present paper is to consider a natural
$q$--analog for (\ref{2.7}). Namely, let us define the following
$q$--analog of the sum (\ref{2.7})
\begin{equation}\sum_{\ld}q^{\frac{1}{2}{\wt\ld}B{\wt\ld}^t}\prod_j\left[
 \begin{array}{c}
    ((E-B){\wt\ld}^t+b^t)_j  \\
    \ld_j
 \end{array}\right]_{q^{\epsilon_j}},
 \end{equation}
where ~~$\epsilon_j=(-1)^{r(j)}$.

Let us recall that $\left[
 \begin{array}{c}
    M \\
    N
 \end{array}\right]_q$ is the Gaussian $q$--binomial coefficient:
 $$\left[
 \begin{array}{c}
    M \\
    N
 \end{array}\right]_q =\left\{
 \begin{array}{ll}
    \ds\frac{(q;q)_M}{(q;q)_N(q;q)_{M-N}}, & \mbox{if $0\le N\le M$},  \\ \\
    0, & \mbox{otherwise}.
 \end{array}\right.
 $$

\begin{rem} {\rm In our previous paper \cite{KL}, see (5.1) and
(5.2), we had considered another $q$--analog of (\ref{2.7}). It
turned out however that the $q$--analog (5.1) from~\cite{KL},
probably, does not possess good combinatorial properties.}
\end{rem}

One of the main results of the present paper is the following:
\medbreak

\begin{theorem}\label{t2.3} Assume that $p_0\ge 1$
be a rational number,  
and  consider a rational function
\begin{equation}\label{2.9}
V_l(q):= V_l^{(p_0)}(q)=\sum_{\ld}
\frac{q^{\frac{1}{2}{\wt\ld}B{\wt\ld}^t}}{\ds\prod_j(q^{\epsilon(j)};q^{\epsilon(j)})_{\ld_j}},
\label{8}
\end{equation}
summation in (\ref{2.9}) is taken over all sequences of
non--negative integer numbers\break $\ld =\{\ld_k\}_{k=1}^{m_{\al
+1}}$ such that
$$l=\sum_{k\ge 1}n_k\ld_k,~~\ld_k\ge 0.
$$
Then we have
\begin{eqnarray}\label{2.11b}
&&\sum_{l\ge 0}q^{\frac{l^2}{p_0}}V_l^{(p_0)}(q)=\\
&&~~~~~~~~~1+\sum_{\scriptsize{\begin{array}{c}k\ge m\ge 0\\
(k,m)\ne
(0,0)\end{array}}}(-1)^{\frac{1+(-1)^{\al}}{2}k+m}q^{(ky_{\al +1}
+my_{\al})(kz_{\al}+mz_{\al -1})+\Delta_{\al}(k,m)}\cdot
\frac{Q_{k,m}^{((-1)^\al)}(q)}{(q;q)_k},\nonumber
\end{eqnarray}
where ~~$$Q_{k,m}^{(\pm 1)}(q):= Q_{k,m}^{(\pm
)}(q)=q^{\frac{1\mp 1}{2}(k+m)} \left[\begin{array}{c}k-1\\
m\end{array}\right]_q
+q^{\frac{1\pm 1}{2}(2k-m)}\left[\begin{array}{c}k-1\\
m-1\end{array}\right]_q,
$$
\vskip 0.2cm
$$\Delta_{\al}(k,m)=\frac{1+(-1)^{\al}}{2}\left(
\begin{array}{c}k-m\\ 2\end{array}\right)
+\frac{1-(-1)^{\al}}{2}\left(\begin{array}{c}m\\
2\end{array}\right).
$$
\end{theorem}
Let us emphasize that polynomials $Q_{k,m}^{(\pm)}(q)$ are the
{\it same} for all rational numbers $p_0\ge 1$.

Taking the sum with respect to the index $m$ in the
RHS(\ref{2.11b}), we obtain
\begin{cor} {\rm (Rogers--Ramanujan--Gordon--Andrews' type identity)}
\begin{equation}\label{2.12}
\sum_{l\ge 0}q^{\frac{l^2}{p_0}}V_l^{(p_0)}(q)=
1+\sum_{k>0}(-1)^{\frac{1-(-1)^{\al}}{2}k}q^{k^2(y_{\al
+1}+y_{\al})(z_{\al} +z_{\al -1})}\cdot
\frac{L_k^{(p_0)}(q)}{(q;q)_k},
\end{equation}
where
\begin{eqnarray*} &&L_k^{(p_0)}(q)=\sum_{m\ge
0}(-1)^mq^{m^2y_{\al}z_{\al -1} -km(y_{\al +1}z_{\al
-1}+2y_{\al}z_{\al -1}+y_{\al}z_{\al})+\Delta_{\al}(k,k-m)}\,
Q_{k,k-m}^{((-1)^{\al})}(q).
\end{eqnarray*}
\end{cor}
  \vskip 0.2cm
A proof of identity (\ref{2.11b}) is a "$q$--version" of that
given in \cite{KL}, Theorem 4.1.

\begin{rem}\label{r2.5} {\rm (Gordon--Andrews' type identity)
Let $p_0\ge 1$ be an integer, then $\al =0$, $y_1=p_0$, $y_0=1$,
$z_0=1$, $z_{-1}=0$, and the RHS of (\ref{2.11b}) takes the
following form
\begin{eqnarray}
&&1+\sum_{\scriptsize{\begin{array}{c}k\ge m\ge 0\\
(k,m)\ne
(0,0)\end{array}}}(-1)^{k+m}q^{(kp_0+m)k+\frac{(k-m)(k-m-1)}{2}}
\cdot\frac{Q_{k,m}^{(+)}(q)}{(q;q)_k}\nn\\
&&~~
=1+\sum_{k>0}(-1)^k\frac{q^{k^2p_0+\frac{k(k-1)}{2}}}{(q;q)_k}
\left\{\sum_{m=0}^k(-1)^mq^{\frac{m(m+1)}{2}}\,
Q_{k,m}^{(+)}(q)\right\}.\label{2.13}
\end{eqnarray}

It is not difficult to see that if $k>0$, the sum in the brackets
(\ref{2.13}) is equal to $(1+q^k)(q;q)_k$. Hence, if $p_0\ge 1$ is
an integer, then we come to the following identity:
\begin{equation}\label{2.14}
\sum_{l\ge 0}q^{\frac{l^2}{p_0}}V_l^{(p_0)}(q)=1+\sum_{k>
0}(-1)^kq^{k^2p_0+\frac{k(k-1)}{2}}(1+q^k).
\end{equation}
Using the Jacobi triple identity, one can rewrite (\ref{2.14}) in
the following forms
\begin{eqnarray}
&&\sum_{l\ge 0}q^{\frac{l^2}{p_0}}V_l^{(p_0)}(q)
=\prod_{n\ge 1}(1-q^{(2p_0+1)n})(1-q^{(2p_0+1)n-p_0-1})
(1-q^{(2p_0+1)n-p_0}),~~~~~~~~\\~~\nn\\
&&\frac{1}{(q;q)_{\infty}}\sum_{l\ge
0} q^{\frac{l^2}{p_0}}V_l^{(p_0)}(q)=\prod_{n\not\equiv
0,p_0,p_0+1(\mod 2p_0+1)}(1-q^n)^{-1}.
\end{eqnarray}

It looks very challenging task to find for any rational number
$p_0\ge 1$ an explicit product formula (Weyl's denominator
identity) for the right hand side of identity (\ref{2.11b}). }
\end{rem}

We consider the identity (\ref{2.12}) as an identity between {\it
bosonic} and {\it fermionic} formulae for the character of "vacuum
representation" of the generalized Kac--Moody algebra
corresponding to the matrix $\Theta^{-1}$.

\section{$XXZ\to XXX$ bijection}
\label{xxzb}
\neweq

In this section we are going to describe a bijection between the
space of states for $XXZ$--model and that of $XXX$--model. Let us
formulate the corresponding combinatorial problem more explicitly.
First of all as it follows from the results of our previous
paper, the combinatorial completeness of Bethe's states for the
$XXZ$ model is equivalent to the following identity
\begin{equation}
\prod_m(2s_m+1)^{N_m}=\sum_{l=0}^{N}Z^{XXZ}(N,s~|~l),
\end{equation}
where $N=\ds\sum_m2s_mN_m$ and the numbers $Z^{XXZ}(N,s~|~l)$ are
given by (\ref{2.7}). On the other hand it follows from the
combinatorial completeness of Bethe's states for $XXX$ model (see
\cite{K1}) that
\begin{equation}
\prod_m(2s_m+1)^{N_m}=\sum_{l=0}^{\frac{1}{2}N}(N-2l+1)Z^{XXX}(N,s~|~l),
\end{equation}
where the number $Z^{XXX}(N,s~|~l)$ stands for the multiplicity of
$\left(\frac{N}{2}-l\right)$-spin irreducible representation of
$sl(2)$ in the tensor product
$$V_{s_1}^{\otimes N_1}\otimes\cdots\otimes V_{s_m}^{\otimes N_m}.
$$

Let us remark that both numbers $Z^{XXZ}(N,s~|~l)$ and
$Z^{XXX}(N,s~|~l)$ admit a combinatorial interpretation in terms
of rigged configurations. The difference between the space of
states of $XXX$ model and that of $XXZ$ model is the availability
of the so--called $1^{-}$--configurations (or $1^-$ string) in
the space of states for the latter model. The presence of
$1^-$--strings in the space of states for $XXZ$--model is a
consequence of broken $sl(2)$--symmetry of the $XXZ$--model. Our
goal in this section is to understand from a combinatorial point
of view how the anisotropy of $XXZ$ model breaks the symmetry of
the $XXX$ chain. More exactly, we suppose to describe a bijection
between $XXZ$--rigged configurations and $XXX$--rigged
configurations. Let us start with recalling a definition of
rigged configurations.

We consider at first the case of $sl(2)$ $XXX$--magnet. Given a
composition $\mu =(\mu_1,\mu_2,\ldots )$ and a natural integer
$l$, by definition a $sl(2)$--configuration of type $(l,\mu )$ is
a partition $\nu\vdash l$ such that all vacancy numbers
\begin{equation}
P_n(\nu ;\mu ):=\sum_k\min (n,\mu_k)-2\sum_{k\le n}\nu'_k
\end{equation}
are nonnegative. Here $\nu'$ denotes the conjugate partition to
that $\nu$. A rigged configuration of type $(l,\mu )$ is a
configuration $\nu$ of type $(l,\mu )$ together with the
collection of integer numbers $\{J_{\alpha}\}_{\alpha
=1}^{m_n(\nu )}$ which satisfy the following inequalities
$$0\le J_1\le J_2\le\cdots\le J_{m_n(\nu )}\le P_n(\nu ;\mu ).
$$
Here $m_n(\nu )$ is equal to the number of parts of the partition
$\nu$ which are equal to $n$. It is clear that the total number
of rigged configurations of type $(l,\mu )$ is equal to the
following number
$$Z(l~|~\mu ):=\sum_{\nu\vdash l}\prod_{n\ge 1}\left (
\begin{array}{c}
    P_n(\nu ;\mu )+m_n(\nu )  \\
    m_n(\nu )
\end{array}\right ).
$$
The following result had been proved in \cite{K1}.

\begin{theorem} Multiplicity of $(N-2l+1)$--dimensional
irreducible representation of $sl(2)$ in the tensor product
$$V_{s_1}^{\otimes N_1}\otimes\cdots\otimes V_{s_m}^{\otimes N_m}
$$
is equal to the number ~~$Z\left(l~|~\underbrace{2s_1,\ldots
,2s_1}_{N_1},\ldots ,\underbrace{2s_m,\ldots ,2s_m}_{N_m}\right)$.
\end{theorem}

\hskip -0.65cm{\bf Example 2} \ One can check that
$$V_1^{\otimes 5}=6V_0+15V_1+15V_2+10V_3+4V_4+V_5.
$$
In our case we have $\mu =(2^5)$. Let us consider $l=5$. It turns
out that there exist three configurations of type $(3,(2^5))$,
namely

\vskip 0.5cm \hskip 3cm \hbox{
\hbox{\Fsquare(0.4cm,{})\Addsquare(0.4cm,{})\Addsquare(0.4cm,{})\Addsquare(0.4cm,{})\Addsquare(0.4cm,{})~0}
\hskip 1cm \normalbaselines\m@th\offinterlineskip
    \vtop{\hbox{\Hfourbox({},{},{},{})~0}
          \vskip-0.4pt
          \hbox{\Fsquare(0.4cm,{})~1}} \hskip 1cm
\normalbaselines\m@th\offinterlineskip
    \vtop{\hbox{\Hthreebox({},{},{})~0}
          \vskip-0.4pt
          \hbox{\Htwobox({},{})~2}}
 }
\vskip 0.5cm

Hence, $Z(3~|~(2^5))=1+2+3=6={\rm Mult}_{V_0}\left( V_1^{\otimes
5}\right)$. \vskip 0.5cm

Now let us give a definition of $sl(2)$--$XXZ$ configuration. We
consider in this Section only the case when the anisotropy
parameter $p_0$ is an integer, $p_0\in{\bf Z}_{\ge 2}$. Under this
assumption the formulae (\ref{2.5}) and (\ref{2.6}) take the
following form:
\begin{eqnarray*}
&&n_j=j, \ \hbox{\rm if} \ 1\le j<p_0, \ v_j=+1;\\
&&n_{p_0}=1, \ v_{p_0}=-1;\\
&&2\Phi_{k,2s}=\frac{2sk}{p_0}-\min (k,2s), \ \hbox{\rm if} \ 1\le k<p_0,
\ 2s+1<p_0;\\
&&2\Phi_{p_0,2s}=\frac{2s}{p_0}, \ \hbox{\rm if} \ 2s+1<p_0;\\
&&b_{kj}=k-j, \ \hbox{\rm if} \ 1\le j\le k<p_0;\\
&&b_{kp_0}=1, \ \hbox{\rm if} \ 1\le k<p_0;\\
&&a_j:=a_j(l~|~\mu )=\sum_m\min
(j,\mu_m)-2l-j\left[\frac{\sum_m\mu_m-2l}{p_0}\right], \ \hbox{\rm if} \ 1\le
j<p_0;\\
&&a_{p_0}(l~|~\mu )=\left[\frac{\sum_m\mu_m-2l}{p_0}\right].
\end{eqnarray*}

\begin{de} A $sl(2)$--$XXZ$--configuration of type $(l,\mu )$ is a pair
$(\ld ,\ld_{p_0})$, where $\ld$ is a composition with all parts strictly
less than $p_0$, 
$\ds\sum_{j<p_0}j\ld_j +\ld_{p_0}=l$, and such that all vacancy numbers
$P_j(\ld ~|~\mu )$ are nonnegative.
\end{de}

Let us recall \cite{KL} that if the anisotropy parameter $p_0\ge
2$ is an integer, then
\begin{eqnarray}
&&P_j(\ld |\mu ):=a_j(l~|~\mu
)+2\ds\sum_{j<k<p_0}(k-j)\ld_k+\ld_{p_0}, \ \hbox{\rm if} \
j<p_0-1;\label{3.4}\\~ \nn \\ &&P_{p_0-1}(\ld ~|~\mu
):=a_{p_0-1}(l~|~\mu)+\ld_{p_0};\nonumber
\\~\nn \\
&&P_{p_0}(\ld ~|~\mu ):=a_{p_0}(l~|~\mu )+\ld_{p_0-1}.\nonumber
\end{eqnarray}

In the sequel we are displaying a configuration
$$(\ld ,\ld_{p_0})=(\ld_1,\ld_2,\ldots ,\ld_{p_0-1},\ld_{p_0})
$$
as the diagram of the following partition
$(1^{\ld_1+\ld_{p_0}},2^{\ld_2},\ldots ,(p_0-1)^{\ld_{p_0-1}})$.

\vskip 0.3cm

\hskip -0.65cm{\bf Example 3} \ Let us consider $p_0=6$,
$s=\frac{3}{2}$, $N=5$, $l=5$. The total number of type
$(5,(3^5))$ $sl(2)$--$XXZ$ configurations is equal to $12$.

\vskip 0.8cm \hbox{\hskip 3.2cm
\hbox{\Vfivebox(\hbox{$\clubsuit$},\hbox{$\clubsuit$},\hbox{$\clubsuit$},
\hbox{$\clubsuit$},\hbox{$\clubsuit$})~0} \hskip 1cm
\hbox{\normalbaselines\m@th\offinterlineskip
    \vtop{\hbox{\Fsquare(0.4cm, )~3}
          \vskip-0.4pt
          \hbox{\Fsquare(0.4cm,\hbox{$\clubsuit$})~0}
          \vskip-0.4pt
          \hbox{\Fsquare(0.4cm,\hbox{$\clubsuit$})}
          \vskip-0.4pt
          \hbox{\Fsquare(0.4cm,\hbox{$\clubsuit$})}
          \vskip-0.4pt
          \hbox{\Fsquare(0.4cm,\hbox{$\clubsuit$})}}}
\hskip 1cm \hbox{\normalbaselines\m@th\offinterlineskip
    \vtop{\hbox{\Fsquare(0.4cm, )~1}
          \vskip-0.4pt
          \hbox{\Fsquare(0.4cm, )}
          \vskip-0.4pt
          \hbox{\Fsquare(0.4cm,\hbox{$\clubsuit$})~0}
          \vskip-0.4pt
          \hbox{\Fsquare(0.4cm,\hbox{$\clubsuit$})}
          \vskip-0.4pt
          \hbox{\Fsquare(0.4cm,\hbox{$\clubsuit$})}}}
\hskip 1cm \normalbaselines\m@th\offinterlineskip
    \vtop{\hbox{\Htwobox({},{})~6}
          \vskip-0.4pt
          \hbox{\Fsquare(0.4cm,\hbox{$\clubsuit$})~0}
          \vskip-0.4pt
          \hbox{\Fsquare(0.4cm,\hbox{$\clubsuit$})}
          \vskip-0.4pt
          \hbox{\Fsquare(0.4cm,\hbox{$\clubsuit$})}}
\hskip 1cm \normalbaselines\m@th\offinterlineskip
    \vtop{\hbox{\Htwobox({},{})~4}
          \vskip-0.4pt
          \hbox{\Fsquare(0.4cm,)~1}
          \vskip-0.4pt
          \hbox{\Fsquare(0.4cm,\hbox{$\clubsuit$})~0}
          \vskip-0.4pt
          \hbox{\Fsquare(0.4cm,\hbox{$\clubsuit$})}}
          }
 \vskip 1cm
\hbox{\hskip 3.2cm  \normalbaselines\m@th\offinterlineskip
    \vtop{\hbox{\Hthreebox({},{},{})~9}
          \vskip-0.4pt
          \hbox{\Fsquare(0.4cm,\hbox{$\clubsuit$})~0}
              \vskip-0.4pt
          \hbox{\Fsquare(0.4cm,\hbox{$\clubsuit$})}}
\hskip 0.7cm
\normalbaselines\m@th\offinterlineskip
    \vtop{\hbox{\Htwobox({},{})~2}
          \vskip-0.4pt
          \hbox{\Htwobox({},{})}
              \vskip-0.4pt
          \hbox{\Fsquare(0.4cm,\hbox{$\clubsuit$})~0}}
\hskip 0.7cm
\normalbaselines\m@th\offinterlineskip
    \vtop{\hbox{\Hthreebox({},{},{})~7}
          \vskip-0.4pt
          \hbox{\Fsquare(0.4cm,)~1}
             \vskip-0.4pt
          \hbox{\Fsquare(0.4cm,\hbox{$\clubsuit$})~0}}
\hskip 0.7cm
\normalbaselines\m@th\offinterlineskip
    \vtop{\hbox{\Hfourbox({},{},{},{})~7}
          \vskip-0.4pt
          \hbox{\Fsquare(0.4cm,\hbox{$\clubsuit$})~0}}
                   }
 \vskip 0.9cm
\hbox{ \hskip 3.2cm \normalbaselines\m@th\offinterlineskip
    \vtop{\hbox{\Hfourbox({},{},{},{})~5}
          \vskip-0.4pt
          \hbox{\Fsquare(0.4cm,)~1}}
\hskip 0.7cm
\normalbaselines\m@th\offinterlineskip
    \vtop{\hbox{\Hthreebox({},{},{})~5}
          \vskip-0.4pt
          \hbox{\Htwobox({},{})~2}}
\hskip 0.7cm \Hfivebox( , , , , )~5 } \vskip 0.8cm

The total number of type $(5,(3^5))$ rigged configurations is equal to
$$Z^{XXZ}(5~|~(3^5))=101=1+4+3+7+10+10+6+16+8+12+18+6.
$$
Here we have used the symbol \ $\clubsuit$ \ to mark the
$1^-$--strings.

Now we are ready to describe a map $\Pi$ from the space of states
for $XXZ$ model to that of $XXX$ one. More exactly we are going to
describe a rule how a $XXZ$--configuration falls to the
$XXX$--pieces. At first we describe this
rule schematically: 
\vskip 0.4cm

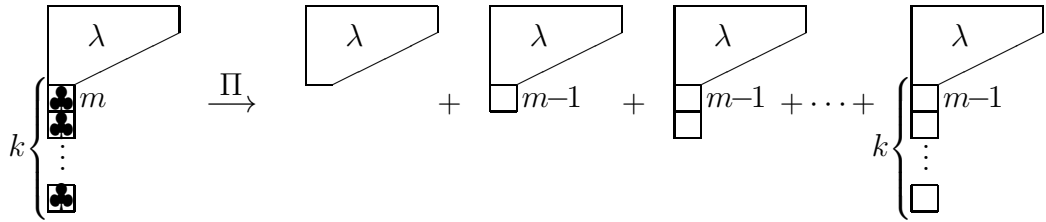
\setlength{\unitlength}{0.35cm}
\begin{picture}(25,11.5)(-1.5,-1.75)
\put(1.2,5){\line(1,0){1}} 
\put(1.2,8){\line(1,0){5}} 
\put(6.2,8){\line(0,-1){1}}
\put(2.2,5){\line(2,1){4}} 
\put(1.2,5){\line(0,1){3}} 
\put(-0.3,2.3){$k\!\left\{\begin{array}{c}~\\~\\~\\~\\ \end{array}\right.$}
\put(1.6,1.8){$\vdots$}
\put(2.4,4.2){$m$}
\put(1.2,3){\line(0,1){2}}
\put(1.2,3){\line(1,0){1}}
\put(2.2,3){\line(0,1){2}}
\put(1.2,4){\line(1,0){1}}
\put(1.25,3.1){$\clubsuit$}
\put(1.25,4.1){$\clubsuit$}
\put(1.25,0.4){$\clubsuit$}
\put(1.2,0.2){\line(1,0){1}}
\put(1.2,0.2){\line(0,1){1}}
\put(1.2,1.2){\line(1,0){1}}
\put(2.2,0.2){\line(0,1){1}}

\put(11,5){\line(1,0){1}} 
\put(11,8){\line(1,0){5}} 
\put(16,8){\line(0,-1){1}}
\put(12,5){\line(2,1){4}} 
\put(11,5){\line(0,1){3}} 

\put(18,5){\line(1,0){1}} 
\put(18,8){\line(1,0){5}} 
\put(23,8){\line(0,-1){1}}
\put(19,5){\line(2,1){4}} 
\put(18,5){\line(0,1){3}} 
\put(18,4){\line(1,0){1}}
\put(18,4){\line(0,1){1}}
\put(19,4){\line(0,1){1}}
\put(19.2,4.2){$m\!\!-\!\!1$}

\put(25,5){\line(1,0){1}} 
\put(25,8){\line(1,0){5}} 
\put(30,8){\line(0,-1){1}}
\put(26,5){\line(2,1){4}} 
\put(25,5){\line(0,1){3}} 
\put(25,3){\line(0,1){2}}
\put(25,3){\line(1,0){1}}
\put(26,3){\line(0,1){2}}
\put(25,4){\line(1,0){1}}
\put(26.2,4.2){$m\!\!-\!\!1$}

\put(34,5){\line(1,0){1}} 
\put(34,8){\line(1,0){5}} 
\put(39,8){\line(0,-1){1}}
\put(35,5){\line(2,1){4}} 
\put(34,5){\line(0,1){3}} 
\put(34,3){\line(0,1){2}}
\put(34,3){\line(1,0){1}}
\put(35,3){\line(0,1){2}}
\put(34,4){\line(1,0){1}}
\put(34,0.2){\line(1,0){1}}
\put(34,0.2){\line(0,1){1}}
\put(34,1.2){\line(1,0){1}}
\put(35,0.2){\line(0,1){1}}
\put(35.2,4.2){$m\!\!-\!\!1$}
\put(32.5,2.3){$k\!\left\{\begin{array}{c}~\\~\\~\\~\\ \end{array}\right.$}
\put(34.4,1.8){$\vdots$}

\put(2.7,6.5){$\lambda$} \put(12.5,6.5){$\lambda$}
\put(19.5,6.5){$\lambda$} \put(26.5,6.5){$\lambda$}
\put(35.5,6.5){$\lambda$}
\put(7.2,4){$\buildrel\displaystyle\Pi\over\longrightarrow$}
\put(16,4){+} \put(23,4){+} \put(29,4){$+\cdots +$}

\end{picture}
\vskip 0.2cm

This decomposition corresponds to the well--known identity
$$\left[
\begin{array}{c}
    m+k\\
    k
\end{array}\right]_q=\sum_{j=0}^kq^j\left[\begin{array}{c}
    m+j-1\\
    j
\end{array}\right]_q.
$$

In what follows we will assume that $p_0>\sum_ms_m$.

\begin{theorem} The map $\Pi$ is well--defined and gives rise to a bijection
between the space of states of $XXZ$--model and that of $XXX$ one.
\end{theorem}

{\it Proof.} Let us start with rewriting the formulae (\ref{3.4})
for the $XXZ$--vacancy numbers in more convenient form, namely,
\begin{eqnarray}
&&P_j^{XXZ}(\wt\nu~|~\mu )=\sum_m\min (j,2s_m)-2\sum_{k\le
j}\nu_k'-j\left[\frac{\sum_m2s_m-2l}{p_0}\right], 
~~~\hbox{\rm if}~~1\le j<p_0-1;\nonumber \\
&& \nonumber \\
&&P_{p_0-1}^{XXZ}(\wt\nu~|~\mu
)=p_0\left\{\frac{\sum_m2s_m-2l}{p_0
}\right\}+\left[\frac{\sum_m2s_m-2l}{p_0}\right] +\ld_{p_0};\label{3.5} \\
~\nn \\
&&P_{p_0}^{XXZ}(\wt\nu~|~\mu )
=\left[\frac{\sum_m2s_m-2l}{p_0}\right] +m_{p_0-1}(\nu ).\nonumber
\end{eqnarray}
Here $\mu =(2s_1,\ldots ,2s_m)$ and $\wt\nu$ is a pair $\wt\nu =(\nu
,\ld_{p_0})$, where $\nu$ is a partition such that $l(\nu )\le p_0-1$,
$|\nu |+\ld_{p_0}=l$. Relationship between $\ld$ from Definition~1 and
$\nu$ is the following
$$m_j(\nu )=\ld_j,~~{\rm i.e.}~~\nu =(1^{\ld_1}2^{\ld_2}\ldots
(p_0-1)^{\ld_{p_0-1}}).
$$
Next, let us consider an integer $l\le\sum_ms_m$ and let
$\nu\vdash l$ be a $XXX$--configuration. Let $\ld_{p_0}$ be an
integer such that $2\sum s_m-2l-p_0<\ld_{p_0}\le\sum_ms_m-l$ and
consider the pair $\wt\nu =(\nu ,\ld_{p_0})$. It is easy to check
that
\begin{eqnarray*}
&&P_j^{XXZ}(\wt\nu~|~\mu )=\sum_m\min (j,2s_m)-2\sum_{k\le
j}\nu_k'=P_j^{XXX}(\nu~|~\mu )\ge 0,
~~~\hbox{\rm if}~~1\le j<p_0-1; \nonumber \\
&& \nonumber \\
&&P_{p_0-1}^{XXZ}(\wt\nu~|~\mu )=\sum_m2s_m-2l+\ld_{p_0}\ge 0;
\\~\nn \\
&&P_{p_0}^{XXZ}(\wt\nu~|~\mu )=\ld_{p_0-1}\ge 0.
\end{eqnarray*}
Thus the pair $\wt\nu =(\nu ,\ld_{p_0})$ is a
$XXZ$--configuration.

Furthermore it follows from our assumptions (namely,
$\sum_ms_m<p_0$,\break $\ld_{p_0}>0$) that $\ld_{p_0-1}=0$ and
both $1^-$--strings and $(p_0-1)$--strings do not give a
contribution to the space of $XXZ$--states. Thus we see that both
$XXX$--configuration $\nu$ and $XXZ$--configuration $\wt\nu =(\nu
,\ld_{p_0})$ define the same number of states. Now, if $\wt\nu
=(\nu ,\mu )$ is a $XXZ$--configuration then $\nu$ is a $XXX$
configuration as well. This is clear because (see (\ref{3.5}))
$$P_j^{XXX}(\nu~|~\mu )\ge P_j^{XXZ}(\wt\nu~|~\mu ),~~1\le j\le p_0-1.
$$
By the similar reasons if $(\wt\nu ,\ld_{p_0})$ is a
$XXZ$--configuration, then for any integer $k$, $0\le
k\le\ld_{p_0}$, the pair $(\wt\nu ,\ld_{p_0}-k)$ is also a
$XXZ$--configuration. It follows from the above considerations
that $\Pi$ is a well--defined map. Furthermore there exists one
to one correspondence between the space of $XXX$--configurations
and that of $XXZ$--configurations, namely,
$$\nu\leftrightarrow\wt\nu =(\nu ,\ld_{p_0}),
$$
where $\ld_{p_0}=[\sum_ms_m-|\nu |]$.

All others $XXZ$--configurations $(\nu ,k)$ with $0\le
k<\sum_ms_m-|\nu |-1$ give a contribution to the space of
descendants for $\nu\leftrightarrow\wt\nu$.

\qed

\vskip 0.8cm


\end{document}